\documentclass{birkjour}
\usepackage{amsmath, amsthm, amscd, amsfonts, amssymb, graphicx, color}
\DeclareMathOperator{\Tr}{Tr}
 \newtheorem{thm}{Theorem}[section]
 \newtheorem{cor}[thm]{Corollary}
 \newtheorem{lem}[thm]{Lemma}
 \newtheorem{prop}[thm]{Proposition}
 \theoremstyle{definition}
 
 \theoremstyle{remark}
 \newtheorem{rem}[thm]{Remark}
 
 \numberwithin{equation}{section}

\begin{document}

%
%
%
%
%
%
%
%
%

\title[$*$-$k$-Ricci-Yamabe soliton]
 {Geometry of $*$-$k$-Ricci-Yamabe soliton and gradient $*$-$k$-Ricci-Yamabe soliton on Kenmotsu manifolds}

\author[S. Dey]{Santu Dey}
\address{Department of Mathematics\\
Bidhan Chandra College\\
Asansol, Burdwan, West Bengal-713304, India.}
\email{santu.mathju@gmail.com}

\author[S. Roy]{Soumendu Roy}
\address{Department of Mathematics\\
Jadavpur University\\
Kolkata-700032, India}
\email{soumendu1103mtma@gmail.com}



\begin{abstract}
The goal of the current paper is to characterize $*$-$k$-Ricci-Yamabe soliton within the framework on Kenmotsu manifolds. Here, we have shown the nature of the soliton and find the scalar curvature when the manifold admitting $*$-$k$-Ricci-Yamabe soliton on Kenmotsu manifold. Next, we have evolved the characterization of the vector field when the manifold satisfies $*$-$k$-Ricci-Yamabe soliton. Also we have embellished some applications of vector field as torse-forming in terms of $*$-$k$-Ricci-Yamabe soliton on Kenmotsu manifold. Then, we have studied gradient $\ast$-$\eta$-Einstein soliton to yield the nature of Riemannian curvature tensor. We have developed an example of $*$-$k$-Ricci-Yamabe soliton on 5-dimensional Kenmotsu manifold to prove our findings.
\end{abstract}

\maketitle
\section{Introduction}
In modern mathematics, the methods of contact geometry make a major contribution. Contact geometry has evolved from the mathematical formalism of classical mechanics. in 1969, S. Tanno \cite{tanno} classified the connected almost contact metric manifolds whose automorphism groups have maximal dimensions as follows:\\
(a) Homogeneous normal contact Riemannian manifolds with constant $\phi$-holomorphic sectional curvature if $k(\xi,X)>0;$\\\\
(b) Global Riemannian product of a line or a circle and a K$\ddot{a}$hlerian manifold with constant holomorphic sectional curvature if $k(\xi,X)=0;$\\\\
(c) A warped product space $\mathbb{R} \times_\lambda \mathbb{C}^n$ if $k(\xi, X) < 0$; where $k(\xi, X)$ denotes the sectional curvature of the plane section containing the characteristic vector field $\xi$ and an arbitrary vector field $X$.\\
In 1972, K. Kenmotsu \cite{kenmotsu} obtained some tensor equations to characterize the manifolds of the third class. Since then the manifolds of the third class were called Kenmotsu manifolds.\\\\
In 1982, R. S. Hamilton \cite{rsham} introduced the concept of Ricci flow, which is an evolution equation for metrics on a Riemannian manifold. The Ricci flow equation is given by:
\begin{equation}\label{1.1}
\frac{\partial g}{\partial t} = -2S,
\end{equation}
on a compact Riemannian manifold $M$ with Riemannian metric $g$.\\
A self-similar solution to the Ricci flow (\cite{rsham}, \cite{topping}) is called a Ricci soliton \cite{rsha} if it moves only by a one parameter family of diffeomorphism and scaling. The Ricci soliton equation is given by \cite{soumendu}:
\begin{equation}\label{1.2}
\pounds_V g + 2S + 2\Lambda g=0,
\end{equation}
where $\pounds_V$ is the Lie derivative in the direction of $V$, $S$ is Ricci tensor, $g$ is Riemannian metric, $V$ is a vector field and $\Lambda$ is a scalar. The Ricci soliton is said to be shrinking, steady and expanding accordingly as $\Lambda$ is negative, zero and positive respectively.\\\\
Recently, Wang-Gomes-Xia \cite{wang} extended the notion of almost Ricci soliton to $k$-almost Ricci soliton. A complete Riemannian manifold $(M^n, g)$ is said to be a $k$-almost Ricci soliton, denoted by $(M^n, g, X, k, \Lambda)$, if there exists smooth vector field $X$ on $M^n$, a soliton function $\Lambda \in C^\infty(M^n)$ and a non-zero real valued function $k$ on $M^n$ such that
\begin{equation}\label{1.2A}
k\pounds_V g + 2S + 2\Lambda g=0.
\end{equation}
The concept of Yamabe flow was first introduced by Hamilton \cite{rsha} to construct Yamabe metrics on compact Riemannian manifolds. On a Riemannian or pseudo-Riemannian manifold $M$, a time-dependent metric $g(\cdot, t)$ is said to evolve by the Yamabe flow if the metric $g$ satisfies the given equation,
\begin{equation}\label{1.3}
  \frac{\partial }{\partial t}g(t)=-rg(t),\hspace{0.5cm} g(0)=g_{0},
\end{equation}
where $r$ is the scalar curvature of the manifold $M$.\\\\
In 2-dimension the Yamabe flow is equivalent to the Ricci flow \cite{rsham} (defined by $\frac{\partial }{\partial t}g(t) = -2S(g(t))$, where $S$ denotes the Ricci tensor). But in dimension $> 2$ the Yamabe and Ricci flows do not agree, since the Yamabe flow preserves the conformal class of the metric but the Ricci flow does not in general.\\
A Yamabe soliton \cite{barbosa} correspond to self-similar solution of the Yamabe flow, is defined on a Riemannian or pseudo-Riemannian manifold $(M, g)$ as:
\begin{equation}\label{1.4}
  \frac{1}{2}\pounds_V g = (r-\Lambda)g,
\end{equation}
where $\pounds_V g$ denotes the Lie derivative of the metric $g$ along the vector field $V$, $r$ is the scalar curvature and $\Lambda$ is a constant. Moreover a Yamabe soliton is said to be expanding, steady, shrinking depending on $\Lambda$ being positive, zero, negative respectively. If $\Lambda$ is a smooth function then \eqref{1.4} is called almost Yamabe soliton \cite{barbosa}. \\

Very recently, Chen \cite{xchen} introduced a new concept, named $k$-almost Yamabe soliton. According to Chen, a Riemannian metric is said to be a $k$-almost Yamabe soliton if there exists a smooth vector field $V$, a $C^\infty$ function $\Lambda$ and a nonzero function $k$ such that
\begin{equation}\label{1.4A}
  \frac{k}{2}\pounds_V g = (r-\Lambda)g,
\end{equation}
holds. If for any smooth function $f$, $V=Df$ then the previous equation is called gradient $k$-almost Yamabe soliton. If $\Lambda$ is constant, then \eqref{1.4A} is called $k$-Yamabe soliton.\\
Since the introduction of Ricci soliton and Yamabe soliton, many authors (\cite{roy}, \cite{roy2}, \cite{roy3}, \cite{roy4}, \cite{roy5}, \cite{ghosh}, \cite{dong}, \cite{joma}, \cite{singh}) have studied these solitons on contact manifolds.\\\\
Recently in 2019, S. G\"uler and M. Crasmareanu \cite{guler} introduced a new geometric flow which is a scalar combination of Ricci and Yamabe flow under the name Ricci-Yamabe map. This flow is also known as Ricci-Yamabe flow of the type $(\alpha, \beta)$.\\
Let $(M^n, g)$ be a Riemannian manifold and $T^{s}_2 (M)$ be the linear space of its symmetric tensor fields of (0, 2)-type and $Riem(M) \subsetneqq T^{s}_2 (M)$  be the infinite space of its Riemannian metrics. In \cite{guler}, the authors have stated the following definition:\\\\
\textbf{Definition 1.1:}\cite{guler} A Riemannian flow on $M$ is a smooth map:
 $$g:I\subseteq \mathbb{R}\rightarrow Riem(M),$$
 where $I$ is a given open interval. We can call it also as time-dependent (or non-stationary) Riemannian metric.\\\\
\textbf{Definition 1.2:}\cite{guler} The map $RY^{(\alpha, \beta, g)}: I \rightarrow T^{s}_2 (M)$ given by:
$$ RY^{(\alpha, \beta, g)}:= \frac{\partial }{\partial t}g(t)+2\alpha S(t)+\beta r(t)g(t),$$
is called the $(\alpha, \beta)$-Ricci-Yamabe map of the Riemannian flow of the Riemannian flow $(M^n, g)$, where $\alpha, \beta$ are some scalas. If $RY^{(\alpha, \beta, g)} \equiv 0$, then $g(\cdot)$ will be called an $(\alpha, \beta)$-Ricci-Yamabe flow.\\
Also in \cite{guler}, the authors characterized that the $(\alpha, \beta)$-Ricci-Yamabe flow is said to be:\\
$\bullet$ Ricci flow \cite{rsham} if $\alpha = 1$, $\beta = 0$.\\
$\bullet$ Yamabe flow \cite{rsha} if $\alpha = 0$, $\beta = 1$.\\
$\bullet$ Einstein flow (\cite{CATINO}, \cite{roy1}) if $\alpha = 1$, $\beta = -1$.\\\\
A soliton to the Ricci-Yamabe flow is called Ricci-Yamabe solitons if it moves only by one parameter group of diffeomorphism and scaling. The metric of the Riemannain manifold $(M^n, g)$, $n >2$ is said to admit $(\alpha, \beta)$-Ricci-Yamabe soliton or simply Ricci-Yamabe soliton (RYS) $(g, V, \Lambda,\alpha, \beta)$ if it satisfies the equation:
\begin{equation}\label{1.5}
  \pounds_V g+2\alpha S+[2\Lambda-\beta r]g=0,
\end{equation}
where $\pounds_V g$ denotes the Lie derivative of the metric $g$ along the vector field $V$, $S$ is the Ricci tensor, $r$ is the scalar curvature and $\Lambda, \alpha,\beta$ are real scalars.\\\\
In the above equation if the vector field $V$ is the gradient of a smooth function $f$ (denoted by $Df$, $D$ denotes the gradient operator) then the equation \eqref{1.5} is called gradient Ricci-Yamabe soliton (GRYS) and it is defined as:
\begin{equation}\label{1.6}
  Hess f+ \alpha S+\Big[\Lambda-\frac{1}{2}\beta r \Big]g=0,
\end{equation}
where $Hess f$ is the Hessian of the smooth function $f$.\\
Moreover the RYS (or GRYS) is said to be expanding, steady or shrinking according as $\Lambda$ is positive, zero, negative respectively. Also if $\Lambda, \alpha, \beta$ become smooth function then \eqref{1.5} and \eqref{1.6} are called almost RYS and almost GRYS respectively.\\

Now using the previous identities \eqref{1.2A}, \eqref{1.4A} and \eqref{1.5}, we define new notation $k$-Ricci-Yamabe soliton($k$-RYS).
A complete Riemannian manifold $(M^n, g)$ is said to be a $k$-almost Ricci-Yamabe soliton, denoted by $(M^n, g, X, k,\\ \Lambda)$, if there exists smooth vector field $X$ on $M^n$, a soliton function $\Lambda \in C^\infty(M^n)$ and a non-zero real valued function $k$ on $M^n$ such that
\begin{equation}\label{1.6a}
k\pounds_V g+2\alpha S+[2\Lambda-\beta r]g=0.
\end{equation}
In the previous equation if the vector field $V$ is the gradient of a smooth function $f$ (denoted by $Df$, $D$ denotes the gradient operator) then the equation \eqref{1.6a} is called gradient almost $k$-Ricci-Yamabe soliton ($k$-GRYS)and it is defined as:
\begin{equation}\label{1.6b}
  kHess f+ \alpha S+\Big[\Lambda-\frac{1}{2}\beta r \Big]g=0.
\end{equation}
Also, $k$-Ricci-Yamabe soliton($k$-RYS) is called\\
$\bullet$ $k$-Ricci soliton (or gradient $k$-Ricci soliton) if $\alpha=1$ and $\beta=0.$\\
$\bullet$ $k$-Yamabe soliton (or gradient $k$-Yamabe soliton) if $\alpha=0$ and $\beta=2.$\\
$\bullet$ $k$-Einstein soliton (or gradient $k$-Einstein soliton) if $\alpha=1$ and $\beta=1.$\\\\
The concept of $*$-Ricci tensor on almost Hermitian manifolds and $*$-Ricci tensor of real hypersurfaces in non-flat complex space were introduced by Tachibana \cite{tachi} and Hamada \cite{hama} respectively where the $*$-Ricci tensor is defined by:
\begin{equation}\label{1.7}
  S^*(X,Y)=\frac{1}{2}(\Tr \{\varphi\circ R(X, \varphi Y)\})
\end{equation}
for all vector fields $X,Y$ on $M^n$, $\varphi$ is a (1, 1)-tensor field and $\Tr$ denotes Trace.\\
If $S^*(X, Y)=\lambda g(X, Y)+\nu \eta(X)\eta(Y)$ for all vector fields $X, Y$ and $\lambda$, $\nu$ are smooth functions, then the manifold is called $*$-$\eta$-Einstein manifold.\\
Further if $\nu=0$ i.e $S^*(X,Y)=\lambda g(X, Y)$ for all vector fields $X,Y$ then the manifold becomes $*$-Einstein.\\
In 2014, Kaimakamis and Panagiotidou \cite{kaipan} introduced the notion of $*$-Ricci soliton which can be defined as:
\begin{equation}\label{1.8}
  \pounds_V g + 2S^* + 2\Lambda g=0
\end{equation}
for all vector fields $X,Y$ on $M^n$ and $\Lambda$ being a constatnt.\\\\
In \cite{venk}, authors have considered $\ast$-Ricci solitons and gradient almost $\ast$-Ricci solitons on Kenmotsu manifolds and obtained some beautiful results. Very recently, Dey et al. \cite{roy, roy*, roy5, roy6, dey3} have studied $\ast$-Ricci solitons and their generalizations in the framework of almost contact geometry.\\
Recently D. Dey \cite{diba} introduced the notion of $*$-Ricci-Yamabe soliton ($*$-RYS) as follows :\\\\
\textbf{Definition 1.3:} A Riemannian or pseudo-Riemannian manifold $(M, g)$ of dimension $n$ is said to admit $*$-Ricci-Yamabe soliton ($*$-RYS) if
\begin{equation}\label{1.8A}
  \pounds_V g+2\alpha S^*+[2\Lambda-\beta r^*]g=0,
\end{equation}
where $r^* = \Tr(S^*)$ is the $*$- scalar curvature and $\Lambda, \alpha, \beta$ are real scalars.\\
The $*$-Ricci-Yamabe soliton ($*$-RYS) is said to be expanding, steady, shrinking depending on $\Lambda$ being positive, zero, negative respectively. If the vector field $V$ is of gradient type i.e $V = grad(f)$, for $f$ is a smooth function on $M$, then the equation \eqref{1.9} is called gradient $*$-Ricci-Yamabe soliton ($*$-GRYS).\\
Using \eqref{1.2A} and \eqref{1.8A}, we introduce $*$-$k$-Ricci-Yamabe soliton ($*$-$k$-RYS) as:\\\\
\textbf{Definition 1.4:} A Riemannian or pseudo-Riemannian manifold $(M, g)$ of dimension $n$ is said to admit $*$-$k$-Ricci-Yamabe soliton ($*$-$k$-RYS) if
\begin{equation}\label{1.9}
  k\pounds_V g+2\alpha S^*+[2\Lambda-\beta r^*]g=0,
\end{equation}
where $\pounds_V g$ denotes the Lie derivative of the metric $g$ along the vector field $V$, $S^*$ is the $*$-Ricci tensor, $r^* = \Tr(S^*)$ is the $*$-scalar curvature and $\Lambda, \alpha, \beta$ are real scalars.\\
If the vector field $V$ is of gradient type i.e., $V = grad(f)$, for $f$ is a smooth function on $M$, then the soliton equation changes to
\begin{equation}\label{A2}
  kHess f+S^\ast+(\lambda-\frac{r^\ast}{2})g=0
\end{equation}
then the equation \eqref{A2} is called gradient $*$-$k$-Ricci-Yamabe soliton ($*$-$k$-GRYS). The $*$-$k$-Ricci-Yamabe soliton or gradient $*$-$k$-Ricci-Yamabe soliton is said to be expanding, steady, shrinking depending on $\Lambda$ being positive, zero, negative respectively. A $*$-$k$-RYS (or $*$-$k$-GRYS) is called an almost $*$-$k$-RYS (or $*$-$k$-GRYS) if $\alpha, \beta$ and $\Lambda$ are smooth functions on $M$.\\\\
The above notation generalizes a large class of solitons. We can also define some solitons in the following way. A $*$-$k$-RYS (or $*$-$k$-GRYS) is called\\\\
$\bullet$ $*$-$k$-Ricci soliton (or gradient $*$-$k$-Ricci soliton) if $\alpha=1$ and $\beta=0.$\\
$\bullet$ $*$-$k$-Yamabe soliton (or gradient $*$-$k$-Yamabe soliton) if $\alpha=0$ and $\beta=2.$\\
$\bullet$ $*$-$k$-Einstein soliton (or gradient $*$-$k$-Einstein soliton) if $\alpha=1$ and $\beta=1.$\\
On the other hand, a nowhere vanishing vector field $\tau$ on a Riemannian or pseudo-Riemannian manifold $(M, g)$ is called torse-forming \cite{yano} if
\begin{equation}\label{1.10}
\nabla_X \tau = \psi X+\omega(X)\tau,
\end{equation}
where $\nabla$ is the Levi-Civita connection of $g$, $\psi$ is a smooth function and $\omega$ is a 1-form.
Moreover, The vector field $\tau$ is called concircular (\cite{chen}, \cite{kyano}) if the 1-form $\omega$ vanishes identically in the equation \eqref{1.10}. The vector field $\tau$ is called concurrent (\cite{Schouten}, \cite{kyano1}) if in \eqref{1.10} the 1-form $\omega$ vanishes identically and the function $\psi = 1$. The vector field $\tau$ is called recurrent if in \eqref{1.10} the function $\psi = 0$. Finally if in \eqref{1.10} $\psi = \omega = 0$, then the vector field $\tau$ is called a parallel vector field.\\
In 2017, Chen \cite{chen1} introduced a new vector field called torqued vector field. If the vector field $\tau$ staisfies \eqref{1.10} with $\omega(\tau) = 0$, then $\tau$ is called torqued vector field. Also in this case, $\psi$ is known as the torqued function and the 1-form $\omega$ is the torqued form of $\tau$.\\\\
Recent years, many authors have been studied Ricci soliton, Ricci-Yamabe soliton, $*$-Ricci-Yamabe soliton and its characterizations on contact geometry. First, Sharma \cite{Sharma} initiated the study of Ricci solitons in contact geometry. D. Dey \cite{diba} studied new type of soliton namely $*$-Ricci-Yamabe soliton on contact geometry. Roy et al. \cite{roy7} have strained conformal Ricci-Yamabe soliton on perfect fluid space time. Siddiqi and Akyol \cite{siddhi} have discussed the notion of $\eta$-Ricci-Yamabe soliton to set up the geometrical structure on Riemannian submersions admitting $\eta$-Ricci-Yamabe soliton with the potential field. Very recently, Yolda\c{s} \cite{yol} measured Kenmotsu metric in terms of $\eta$-Ricci-Yamabe soliton. Next, Chen \cite{chen2} considered a real hypersurface of a non-flat complex space form which admits a $*$-Ricci soliton whose potential vector field belongs to the principal curvature space and the holomorphic distribution.  Recently, Wang \cite{Wang1} proved that if the metric of a Kenmotsu 3-manifold represents a $*$-Ricci soliton, then the manifold is locally isometric to the hyperbolic space $\mathbb{H}^{3}(-1)$.\\\\
In the ongoing paper, we will discuss about $*$-$k$-Ricci-Yamabe soliton and gradient $*$-$k$-Ricci-Yamabe soliton on Kenmotsu manifold.\\\\
The outline of the article goes as follows:\\
In section 2, after a brief introduction, we have discussed some preliminaries of kenmotsu manifold. In section 3, we have studied $*$-$k$-Ricci-Yamabe soliton admitting Kenmotsu manifold and obtained the nature of soliton, Laplacian of the smooth function. We have also proved that the manifold is $\eta$-Einstein when the manifold satisfies $*$-$k$-Ricci-Yamabe soliton  and the vector field is conformal Killing. Next section, we have demonstrated some applications of torse-forming potential vector field admitting $*$-$k$-Ricci-Yamabe soliton on Kenmotsu manifold. Section 5 deals with gradient $*$-$k$-Ricci-Yamabe soliton and obtain the curvature tensor. In section 6, we have constructed an example to illustrate the existence of $*$-$k$-Ricci-Yamabe soliton on 5-dimensional Kenmotsu manifold.
\section{Preliminaries}
Let $M$ be a (2$n$+1) dimensional connected almost contact metric manifold with an almost contact metric structure $(\phi, \xi, \eta, g)$ where $\phi$ is a $(1, 1)$ tensor field, $\xi$ is a vector field, $\eta$ is a 1-form  and $g$ is the compatible Riemannian metric such that\\
\begin{equation}\label{2.1}
\phi^2(X) = -X + \eta(X)\xi, \eta(\xi) = 1, \eta \circ \phi = 0, \phi \xi = 0,
\end{equation}\\
\begin{equation}\label{2.2}
g(\phi X, \phi Y) = g(X,Y) - \eta(X)\eta(Y),
\end{equation}\\
\begin{equation}\label{2.3}
g(X, \phi Y) = -g(\phi X,Y),
\end{equation}\\
\begin{equation}\label{2.4}
g(X, \xi) = \eta(X),
\end{equation}\\
for all vector fields $X, Y \in \chi(M).$\\
An almost contact metric manifold is said to be a Kenmotsu manifold \cite{kenmotsu} if
\begin{equation}\label{2.5}
  (\nabla_X \phi)Y= -g(X,\phi Y)\xi-\eta(Y)\phi X,
\end{equation}
\begin{equation}\label{2.6}
  \nabla_X \xi=X-\eta(X)\xi,
\end{equation}
where $\nabla$ denotes the Riemannian connection of $g$.\\
In a Kenmotsu manifold the following relations hold (\cite{bagewadi}, \cite{roy4}):
\begin{equation}\label{2.7}
  \eta(R(X,Y)Z)= g(X,Z)\eta(Y)-g(Y,Z)\eta(X),
\end{equation}
\begin{equation}\label{2.8}
  R(X,Y)\xi=\eta(X)Y-\eta(Y)X,
\end{equation}
\begin{equation}\label{2.9}
  R(X,\xi)Y=g(X,Y)\xi-\eta(Y)X,
\end{equation}
where $R$ is the Riemannian curvature tensor.\\
\begin{equation}\label{2.10}
  S(X,\xi)=-2n\eta(X),
\end{equation}
\begin{equation}\label{2.11}
  S(\phi X,\phi Y)= S(X,Y)+2n\eta(X)\eta(Y),
\end{equation}
\begin{equation}\label{2.12}
  (\nabla_X \eta)Y= g(X,Y)-\eta(X)\eta(Y)
\end{equation}
for all vector fields $X, Y, Z \in \chi(M).$\\
Now, we know
\begin{equation}\label{2.13}
  (\pounds_\xi g)(X,Y)=g(\nabla_X \xi, Y)+g(X, \nabla_Y \xi),
\end{equation}
for all vector fields $X, Y, \in \chi(M).$\\
Then using \eqref{2.6} and \eqref{2.13}, we obtain
\begin{equation}\label{2.14}
  (\pounds_\xi g)(X,Y)=2[g(X,Y)-\eta(X)\eta(Y)].
\end{equation}
\begin{prop}
\cite{venk}  On a $(2n+1)$-dimensional Kenmotsu manifold, the $*$-Ricci tensor is given by
\begin{equation}\label{2.15}
  S^*(X,Y) = S(X,Y) + (2n - 1)g(X,Y) + \eta(X) \eta(Y).
\end{equation}
\end{prop}
Also, we take $X = e_i$, $Y = e_i$ in the above equation, where $e_i$'s are a local orthonormal frame and summing over $i=1,2,....,(2n+1)$ to achieve
\begin{equation}\label{2.15new}
  r^*=r+4n^2,
\end{equation}
where $r^*$ is the $*$- scalar curvature of $M$.
\section{Kenmotsu metric as $*$-$k$-Ricci-Yamabe soliton}
Let $M$ be a (2$n$+1) dimensional Kenmotsu manifold. Consider $V = \xi$ in the equation of $*$-$k$-Ricci-Yamabe soliton \eqref{1.9} on $M$, we obtain:
\begin{equation}\label{3.1}
k(\pounds_\xi g)(X,Y)+2\alpha S^*(X,Y)+[2\Lambda-\beta r^*]g(X,Y)=0
\end{equation}
for all vector fields $X, Y, \in \chi(M).$\\
From \eqref{2.14} and \eqref{2.15}, the above equation becomes
\begin{equation}\label{3.2}
  \alpha S(X,Y)+[\Lambda+\alpha(2n - k)+1-\frac{\beta r^*}{2}]g(X,Y)+[\alpha-k]\eta(X)\eta(Y)=0.
\end{equation}
Now, we plug $Y=\xi$ in the above equation and from the identities \eqref{2.1}, \eqref{2.10} to yield
\begin{equation}\label{3.3}
  [\Lambda-\frac{\beta r^*}{2}]\eta(X)=0.
\end{equation}
Since $\eta(X) \neq 0$, the above equation takes the form
\begin{equation}\label{3.4}
  \Lambda=\frac{\beta r^*}{2}.
  \end{equation}
Using \eqref{2.15new}, we acquire
\begin{equation}\label{3.4new}
  \Lambda=\frac{\beta (r+4n^2)}{2}.
\end{equation}
  This leads to the following:
\begin{thm}
If the metric $g$ of a (2$n$+1) dimensional Kenmotsu manifold satisfies the $*$-$k$-Ricci-Yamabe soliton $(g, \xi, \Lambda, \alpha, \beta)$, where $\xi$ is the reeb vector field, then the soliton is expanding, steady, shrinking according as $\beta (r+4n^2) \gtreqqless 0$.
\end{thm}
Also we have, if the manifold $M$ becomes flat i.e., $r = 0$ then \eqref{3.4new} becomes, $\Lambda=2\beta n^2$.\\
So we can state
\begin{cor}
If the metric $g$ of a (2$n$+1) dimensional Kenmotsu manifold, which is flat, satisfies the $*$-$k$-Ricci-Yamabe soliton $(g, \xi, \Lambda, \alpha, \beta)$, where $\xi$ is the reeb vector field, then the soliton is expanding, steady, shrinking according as $\beta \gtreqqless 0$.
\end{cor}
Now consider a $*$-$k$-Ricci-Yamabe soliton $(g, V, \Lambda, \alpha, \beta)$ on $M$ as:
\begin{equation}\label{3.5}
  k(\pounds_V g)(X,Y)+2\alpha S^*(X,Y)+[2\Lambda-\beta r^*]g(X,Y)=0
\end{equation}
for all vector fields $X, Y, \in \chi(M).$\\
Taking $X = e_i$, $Y = e_i$, in the above equation, where $e_i$'s are a local orthonormal frame and summing over $i=1,2,....,(2n+1)$ and using
\eqref{2.15new}, we get
\begin{equation}\label{3.6}
div V+\frac{(r+4n^2)}{k}\Big[\alpha-\frac{\beta(2n+1)}{2}\Big]+\frac{\Lambda}{k}(2n+1)=0.
\end{equation}
If we take the vector field $V$ is of gradient type i.e $V = grad(f)$, for $f$ is a smooth function on $M$, then the equation \eqref{3.6} becomes
\begin{equation}\label{3.7}
\Delta(f)=-\frac{(r+4n^2)}{k}\Big[\alpha-\frac{\beta(2n+1)}{2}\Big]-\frac{\Lambda}{k}(2n+1),
\end{equation}
where $\Delta(f)$ is the Laplacian equation satisfied by $f$.\\
So, we can state the following theorem:
\begin{thm}\label{3.3}
If the metric $g$ of a (2$n$+1) dimensional Kenmotsu manifold satisfies the $*$-$k$-Ricci-Yamabe soliton $(g, V, \Lambda, \alpha, \beta)$, where $V$ is the gradient of a smooth function $f$, then the Laplacian equation satisfied by $f$ is
$$\Delta(f)=-\frac{(r+4n^2)}{k}\Big[\alpha-\frac{\beta(2n+1)}{2}\Big]-\frac{\Lambda}{k}(2n+1).$$
\end{thm}
Now if $\alpha=1, \beta=0$, \eqref{1.9} reduces to $*$-$k$-Ricci soliton and \eqref{3.7} takes the form $\Delta(f)=-\frac{(r+4n^2)}{k}-\frac{\Lambda}{k}(2n+1)$.\\
If $\alpha=0, \beta=2$, \eqref{1.9} reduces to $*$-$k$-Yamabe soliton and \eqref{3.7} takes the form, $\Delta(f)=\frac{1}{k}[r+4n^2- \Lambda](2n+1)$.\\
Moreover if $\alpha=\beta=1$, \eqref{1.9} reduces to $*$-$k$-Einstein soliton and \eqref{3.7} takes the form
$\Delta(f)=-\frac{(r+4n^2)}{k}\Big[1-\frac{(2n+1)}{2}\Big]-\frac{\Lambda}{k}(2n+1)$.\\
Then we have
\begin{rem}
\textbf{Case-I:} If the metric $g$ of a (2$n$+1) dimensional Kenmotsu manifold satisfies the $*$-$k$-Ricci soliton $(g, V, \Lambda)$, where $V$ is gradient of a smooth function $f$, then the Laplacian equation satisfied by $f$ is
$$\Delta(f)=-\frac{(r+4n^2)}{k}-\frac{\Lambda}{k}(2n+1).$$\\
\textbf{Case-II:} If the metric $g$ of a (2$n$+1) dimensional Kenmotsu manifold satisfies the $*$-$k$-Yamabe soliton $(g, V, \Lambda)$, where $V$ is the gradient of a smooth function $f$, then the Laplacian equation satisfied by $f$ is
$$\Delta(f)=\frac{1}{k}[r+4n^2- \Lambda](2n+1).$$\\
\textbf{Case-III:} If the metric $g$ of a (2$n$+1) dimensional Kenmotsu manifold satisfies $*$-$k$-Einstein soliton $(g, V, \Lambda)$,  where $V$ is gradient of a smooth function $f$, then the Laplacian equation satisfied by $f$ is
$$\Delta(f)=-\frac{(r+4n^2)}{k}\Big[1-\frac{(2n+1)}{2}\Big]-\frac{\Lambda}{k}(2n+1).$$
\end{rem}
A vector field $V$ is said to be a conformal Killing vector field iff the following relation holds:
\begin{equation}\label{3.9}
(\pounds_V g)(X,Y)=2\Omega g(X,Y),
\end{equation}
where $\Omega$ is some function of the co-ordinates(conformal scalar).\\
Moreover if $\Omega$ is not constant the conformal Killing vector field $V$ is said to be proper. Also when $\Omega$ is constant, $V$ is called homothetic vector field and when the constant $\Omega$ becomes non zero, $V$ is said to be proper homothetic vector field. If $\Omega = 0$ in the
above equation, then $V$ is called Killing vector field.\\
Let $(g,V,\Lambda, \alpha, \beta)$ be a $*$-Ricci-Yamabe soliton on a (2$n$+1) dimensional Kenmotsu manifold $M$, where $V$ is a conformal Killing vector field. Then from \eqref{1.9}, \eqref{2.15} and \eqref{3.9}, we obtain
\begin{equation}\label{3.10}
  \alpha S(X,Y)=-\Big[\alpha (2n-1)+\Lambda+k\Omega-\frac{\beta r^*}{2}\Big]g(X,Y)-\alpha\eta(X)\eta(Y),
\end{equation}
which leads to the fact that the manifold is $\eta$-Einstein, provided $\alpha \neq 0$.\\
Thus, we have the following theorem:
\begin{thm}
If the metric $g$ of a (2$n$+1) dimensional Kenmotsu manifold satisfies the $*$-$k$-Ricci-Yamabe soliton $(g, V, \Lambda, \alpha, \beta)$, where $V$ is a conformal Killing vector field, then the manifold becomes $\eta$-Einstein, provided $\alpha \neq 0$.
\end{thm}
Taking $Y = \xi$ in the above equation and using \eqref{2.1}, \eqref{2.10}, we achieve
\begin{equation}\label{3.11}
 \Big[ 2\alpha n-\alpha (2n-1)-\Lambda-k\Omega+\frac{\beta r^*}{2}-\alpha \Big]\eta(X)=0.
\end{equation}
Since $\eta(X) \neq 0$, we obtain
\begin{equation}\label{3.12}
  \Omega=\frac{1}{k}\Big[\frac{\beta r^*}{2}-\Lambda\Big].
\end{equation}
Then using \eqref{2.15new}, the above equation becomes
\begin{equation}\label{3.12new}
  \Omega=\frac{1}{k}\Big[\frac{\beta (r+4n^2)}{2}-\Lambda\Big].
\end{equation}
Hence, we can state
\begin{thm}
Let the metric $g$ of a (2$n$+1) dimensional Kenmotsu manifold satisfy the $*$-$k$-Ricci-Yamabe soliton $(g, V, \Lambda, \alpha, \beta)$, where $V$ is
a conformal Killing vector field. Then $V$ is\\
(i)proper vector field if $\frac{1}{k}\Big[\frac{\beta (r+4n^2)}{2}-\Lambda\Big]$ is not constant.\\\\
(ii)homothetic vector field if $\frac{1}{k}\Big[\frac{\beta (r+4n^2)}{2}-\Lambda\Big]$ is constant.\\\\
(iii) proper homothetic vector field if $\frac{1}{k}\Big[\frac{\beta (r+4n^2)}{2}-\Lambda\Big]$ is non-zero constant.\\\\
(iv) Killing vector field if $\Lambda=\frac{\beta (r+4n^2)}{2}$.
\end{thm}
Now, a Kenmotsu manifold $(M^{2n+1}, g)$ is said to $\eta$-Einstein if its Ricci tensor $S$ of type (0, 2) is of the form
\begin{eqnarray}\label{4.6}
S=ag+b\eta\otimes\eta,
\end{eqnarray}
where $a$ and $b$ are smooth function on $(M^{2n+1}, g).$\\
Taking contraction on \eqref{4.6}, we get
\begin{eqnarray}\label{4.7}
r=a(2n+1)+b.
\end{eqnarray}
Now using the identities \eqref{1.5}, \eqref{2.15} and \eqref{2.15new}, we acquire
\begin{eqnarray}\label{4.8}
k[g(\nabla_{X}V, Y)&+&g(X, \nabla_{V}Y)]+2\alpha S(X, Y)+[2\Lambda-\beta(r+4n^2)\nonumber\\
&+&2\alpha(2n-1)]g(X, Y)+2\alpha \eta(X)\eta(Y)=0.
\end{eqnarray}
We insert the identities \eqref{4.6} and \eqref{4.7} into the previous equation to yield
\begin{eqnarray}\label{4.9}
k[g(\nabla_{X}V, Y)&+&g(X, \nabla_{V}Y)]+(2a\alpha+2\alpha(2n-1)+2\Lambda\nonumber\\
&-&\beta(r+4n^2)]g(X, Y)+2\alpha(b+1)\eta(X)\eta(Y)=0.
\end{eqnarray}
Now, we plug $X=Y=\xi$ into \eqref{4.9} to achieve
\begin{eqnarray}\label{4.10}
2kg(\nabla_{\xi}V, \xi)=(2a\alpha+2\alpha(2n-1)+2\Lambda-\beta(r+4n^2)+2\alpha(b+1).
\end{eqnarray}
Putting $V=\xi$ into identity \eqref{4.10}, we obtain
\begin{eqnarray}\label{4.11}
kg(\nabla_{\xi}\xi, \xi)=\Big[a\alpha+2n\alpha+\Lambda-\frac{\beta}{2}(r+4n^2)+\alpha b\Big].
\end{eqnarray}
As it is well known that $$g(\nabla_{\xi}\xi, \xi)=0$$ for all vector field on $M.$ It follows that
$$\Lambda=-a\alpha-2n\alpha+\frac{\beta}{2}\Big(r+4n^2\Big)-b\alpha .$$
So, we have the following theorem:
\begin{thm}
If $(M, \phi, g, \xi, \lambda, a, b)$ is $\ast$-$k$-Ricci-Yamabe soliton on an $\eta$-Einstein Kenmotsu manifold, then $\Lambda=-a\alpha-2n\alpha+\frac{\beta}{2}\Big(r+4n^2\Big)-b\alpha.$
\end{thm}

\section{Application of torse forming vector field on Kenmotsu manifold admitting $*$-$k$-Ricci-Yamabe soliton}
Let $(g, \tau, \Lambda, \alpha, \beta)$ be a $*$-$k$-Ricci-Yamabe soliton on a (2$n$+1) dimensional Kenmotsu manifold $M$, where $\tau$ is a
torse-forming vector field. Then from \eqref{1.9}, \eqref{2.15} and \eqref{2.15new}, we have
\begin{multline}\label{4.1}
  k(\pounds_\tau g)(X,Y)+2\alpha[S(X,Y) + (2n - 1)g(X,Y) + \eta(X)\eta(Y)]\\
   +[2\Lambda-\beta (r+4n^2)]g(X,Y)=0
\end{multline}
where $\pounds_\tau g$ denotes the Lie derivative of the metric $g$ along the vector field $\tau$.\\
Now using\eqref{1.10}, we obtain
\begin{eqnarray}\label{4.2}
  (\pounds_\tau g)(X,Y) &=& g(\nabla_X \tau,Y)+g(X,\nabla_Y \tau) \nonumber\\
                        &=& 2\psi g(X,Y)+\omega(X)g(\tau,Y)+\omega(Y)g(\tau,X)
\end{eqnarray}
for all $X,Y \in M$.\\
Then from \eqref{4.2} and \eqref{4.1}, we get
\begin{multline}\label{4.3}
  \Big[\frac{\beta (r+4n^2)}{2}-\Lambda-k\psi-\alpha(2n-1)\Big]g(X,Y)-\alpha S(X,Y)-\alpha\eta(X)\eta(Y)\\
  =\frac{k}{2}\Big[\omega(X)g(\tau,Y)+\omega(Y)g(\tau,X)\Big].
\end{multline}
Now, we take contraction of \eqref{4.3} over $X$ and $Y$ to acquire
\begin{equation}\label{4.4}
  \Big[\frac{\beta (r+4n^2)}{2}-\Lambda-k\psi-\alpha(2n-1)\Big](2n+1)-\alpha r-\alpha=k\omega(\tau),
\end{equation}
which leads to
\begin{equation}\label{4.5}
  \Lambda=\frac{\beta (r+4n^2)}{2}-k\psi-\alpha(2n-1)-\frac{\alpha r+\alpha+k\omega(\tau)}{(2n+1)}.
\end{equation}
So, we can state the following theorem:
\begin{thm}
 If the metric $g$ of a (2$n$+1) dimensional Kenmotsu manifold satisfies the $*$-$k$-Ricci-Yamabe soliton $(g, \tau, \Lambda, \alpha, \beta)$, where
 $\tau$ is a torse-forming vector field, then $\Lambda=\frac{\beta (r+4n^2)}{2}-k\psi-\alpha(2n-1)-\frac{\alpha r+\alpha+k\omega(\tau)}{(2n+1)}$ and
 the soliton is expanding, steady, shrinking according as $\frac{\beta (r+4n^2)}{2}-k\psi-\alpha(2n-1)-\frac{\alpha r+\alpha+k\omega(\tau)}{(2n+1)}
 \gtreqqless 0$.
 \end{thm}
 Now in \eqref{4.5}, if the 1-form $\omega$ vanishes identically then $\Lambda=\frac{\beta (r+4n^2)}{2}-k\psi-\alpha(2n-1)-\frac{\alpha
 r+\alpha}{(2n+1)}$.\\\\
 If the 1-form $\omega$ vanishes identically and the function $\psi = 1$ in \eqref{4.5}, then $\Lambda=\frac{\beta
 (r+4n^2)}{2}-k-\alpha(2n-1)-\frac{\alpha r+\alpha}{(2n+1)}$.\\\\
 In \eqref{4.5}, if the function $\psi = 0$, then $\Lambda=\frac{\beta (r+4n^2)}{2}-\alpha(2n-1)-\frac{\alpha r+\alpha+k\omega(\tau)}{(2n+1)}$.\\\\
 If $\psi = \omega = 0$ in \eqref{4.5}, then $\Lambda=\frac{\beta (r+4n^2)}{2}-\alpha(2n-1)-\frac{\alpha r+\alpha}{(2n+1)}$.\\\\
 Finally in \eqref{4.5}, if $\omega(\tau) = 0$, then $\Lambda=\frac{\beta (r+4n^2)}{2}-k\psi-\alpha(2n-1)-\frac{\alpha r+\alpha}{(2n+1)}$.\\\\
Then we have
\begin{cor}
Let the metric $g$ of a (2$n$+1) dimensional Kenmotsu manifold satisfy the $*$-$k$-Ricci-Yamabe soliton $(g, \tau, \Lambda, \alpha, \beta)$, where
$\tau$ is a torse-forming vector field, then if $\tau$ is\\
(i) concircular, then $\Lambda=\frac{\beta (r+4n^2)}{2}-k\psi-\alpha(2n-1)-\frac{\alpha r+\alpha}{(2n+1)}$ and the soliton is expanding, steady,
shrinking according as $\frac{\beta (r+4n^2)}{2}-k\psi-\alpha(2n-1)-\frac{\alpha r+\alpha}{(2n+1)} \gtreqqless 0.$\\\\
(ii) concurrent, then $\Lambda=\frac{\beta (r+4n^2)}{2}-1-\alpha(2n-1)-\frac{\alpha r+\alpha}{(2n+1)}$ and the soliton is expanding, steady,
shrinking according as $\Lambda=\frac{\beta (r+4n^2)}{2}-1-\alpha(2n-1)-\frac{\alpha r+\alpha}{(2n+1)} \gtreqqless 0.$\\\\
(iii) recurrent, then $\Lambda=\frac{\beta (r+4n^2)}{2}-\alpha(2n-1)-\frac{\alpha r+\alpha+k\omega(\tau)}{(2n+1)}$ and the soliton is expanding,
steady, shrinking according as $\Lambda=\frac{\beta (r+4n^2)}{2}-\alpha(2n-1)-\frac{\alpha r+\alpha+k\omega(\tau)}{(2n+1)} \gtreqqless 0.$\\\\
(iv) parallel, then $\Lambda=\frac{\beta (r+4n^2)}{2}-\alpha(2n-1)-\frac{\alpha r+\alpha}{(2n+1)}$ and the soliton is expanding, steady, shrinking
according as $\Lambda=\frac{\beta (r+4n^2)}{2}-\alpha(2n-1)-\frac{\alpha r+\alpha}{(2n+1)} \gtreqqless 0.$\\\\
(v) torqued, then $\Lambda=\frac{\beta (r+4n^2)}{2}-k\psi-\alpha(2n-1)-\frac{\alpha r+\alpha}{(2n+1)}$ and the soliton is expanding, steady,
shrinking according as $\Lambda=\frac{\beta (r+4n^2)}{2}-k\psi-\alpha(2n-1)-\frac{\alpha r+\alpha}{(2n+1)} \gtreqqless 0.$
\end{cor}

\section{Gradient $*$-$k$-Ricci-Yamabe soliton on Kenmotsu manifolds}
This section is devoted to the study of Kenmotsu manifolds admitting gradient $*$-$k$-Ricci-Yamabe soliton and we try to characterize the potential vector field of the soliton. First we proof a lemma on Kenmotsu manifold.\\
Now from \eqref{2.15}, we can write
\begin{eqnarray}\label{5.a}
Q^{*}X=-X+\eta(X)\xi
\end{eqnarray}
Differentiating the above equation covariantly with respect to $Y$, we get
\begin{eqnarray}\label{5.b}
\nabla_{Y}Q^{*}X&=&-\nabla_{Y}X+[g(X, Y)-\eta(X)\eta(Y)+\eta(\nabla_{Y}X)-\eta(X)\eta(Y)]\xi\nonumber\\
&+&\nabla_{X}Y.
\end{eqnarray}
Now, we use the above two identities \eqref{5.a} and \eqref{5.b} to yield
\begin{eqnarray}\label{5.c}
(\nabla_{Y}Q^{*})X&=&\nabla_{Y}Q^{*}X-Q^{*}(\nabla_{Y}X)\nonumber\\
&=&g(X, Y)\xi-2\eta(X)\eta(Y)\xi+\eta(X)Y.
\end{eqnarray}
We plug $X=\xi$ into identity \eqref{5.c} to achieve
\begin{eqnarray}\label{5.d}
(\nabla_{Y}Q^{*})\xi=Y-\eta(Y)\xi=\nabla_{Y}\xi.
\end{eqnarray}
Now, we insert $Y=\xi$ into \eqref{5.c} to find
\begin{eqnarray}\label{5.e}
(\nabla_{\xi}Q^{*})Y=0.
\end{eqnarray}
From the identities \eqref{5.d} and \eqref{5.f}, we obtain
\begin{eqnarray}\label{5.f}
(\nabla_{Y}Q^{*})\xi-(\nabla_{\xi}Q^{*})Y=\nabla_{Y}\xi.
\end{eqnarray}
Now we have the following lemma.
\begin{lem}
For a (2n+1)-dimensional Kenmotsu manifold, the following relation holds
$$(\nabla_{Y}Q^{*})\xi-(\nabla_{\xi}Q^{*})Y=\nabla_{Y}\xi.$$
\end{lem}
Since the metric is gradient $*$-$k$-Ricci-Yamabe soliton, so using \eqref{A2}, \eqref{2.15} and \eqref{2.15new}, we can write
\begin{equation}\label{5.1}
k Hess f(X, Y)=-\alpha S(X, Y)-\Big[\Lambda+\frac{\beta(r+4n^2)}{2}+(2n-1)\Big]g(X, Y)- \alpha \eta(X)\eta(Y)
\end{equation}
for all $X, Y \in M$.\\
Now the foregoing equation can be rewritten as
\begin{eqnarray}\label{5.2}
k\nabla_{X}Df=-\alpha QX-\Big[\Lambda+\frac{\beta(r+4n^2)}{2}+(2n-1)\Big]X-\alpha \eta(X)\xi.
\end{eqnarray}
Covariantly diffrentiating the previous equation along an arbitrary vector field $Y$ and using \eqref{2.6}, we achieve
\begin{eqnarray}\label{5.3}
k\nabla_{Y}\nabla_{X}Df+(Yk)\nabla_{X}Df&=&\nabla_{Y}QX-\Big[\lambda+\frac{r+4n^2}{2}-(2n-1)\Big]\nabla_{Y}X\nonumber\\
&-&\frac{Y(r)}{2}X-(\mu-1)[\nabla_{Y}\eta(X)\xi\nonumber\\
&+&(Y-\eta(Y)\xi)\eta(X)].
\end{eqnarray}
Now, we replace $X$ and $Y$ into the identity \eqref{5.3} to yield
\begin{eqnarray}\label{5.4}
k\nabla_{X}\nabla_{Y}Df+(Xk)\nabla_{Y}Df&=&\nabla_{X}QY-\Big[\lambda+\frac{r+4n^2}{2}-(2n-1)\Big]\nabla_{X}Y\nonumber\\
&-&\frac{X(r)}{2}Y-(\mu-1)[\nabla_{X}\eta(Y)\xi\nonumber\\
&+&(X-\eta(X)\xi)\eta(Y)].
\end{eqnarray}
Also in view of \eqref{5.2}, we acquire
\begin{eqnarray}\label{5.5}
k\nabla_{[X, Y]}Df&=&Q(\nabla_{X}Y-\nabla_{Y}X)-\Big[\lambda+\frac{r+4n^2}{2}-(2n-1)\Big](\nabla_{X}Y\nonumber\\
&-&\nabla_{Y}X)-(\mu-1)\eta(\nabla_{X}Y-\nabla_{Y}X)\xi.
\end{eqnarray}
Now, we plug the values of \eqref{5.3}, \eqref{5.4} and \eqref{5.5} into the very well known Riemannian curvature formula
$$R(X, Y)Z=\nabla_{X}\nabla_{Y}Z-\nabla_{Y}\nabla_{X}Z-\nabla_{[X, Y]}Z,$$
to achieve
\begin{eqnarray}\label{5.6}
kR(X, Y)Df&=&\alpha[(\nabla_{Y}Q)X-(\nabla_{X}Q)Y]+\frac{\beta}{2}[X(r)Y-Y(r)X]\nonumber\\
&+&\frac{\alpha}{k}[(Xk)QY-(Yk)QX]+\frac{1}{k}\Big[\Lambda-\frac{\beta(r+4n^2)}{2}\nonumber\\
&+&(2n-1)\Big][(Xk)Y-(Yk)X]+\frac{\alpha}{k}[(Xk)\eta(Y)\nonumber\\
&-&(Yk)\eta(X)]\xi+\alpha[Y\eta(X)-X\eta(Y)].
\end{eqnarray}
So, we have the following lemma:
\begin{lem}
If $(g, V, \alpha, \beta, \Lambda)$ is a gradient $\ast$-$k$-Ricci-Yamabe soliton on a (2n+1)-dimensional Kenmotsu manifold $(M, g, \phi, \xi, \eta)$, then the Riemannian curvature tensor $R$ satisfies \eqref{5.6}.
\end{lem}
\section{Example of a 5-dimensional Kenmotsu manifold admitting $*$-$k$-Ricci-Yamabe soliton}
Let us consider the set $M=\{(x, y, z, u, v)\in\Bbb{R}^5\}$ as our manifold where $(x, y, z, u, v)$ are the standard coordinates in $\Bbb{R}^5$. The vector fields defined below:
  \begin{align*}
    e_1 &= e^{-v}\frac{\partial}{\partial x}, & e_2 &= e^{-v}\frac{\partial}{\partial y}, & e_3 &= e^{-v}\frac{\partial}{\partial z}, & e_4 &= e^{-v}\frac{\partial}{\partial u}, & e_5 &= \frac{\partial}{\partial v}
  \end{align*}
  are linearly independent at each point of $M$. We define the metric $g$ as
  \[
  g(e_i,e_j)=\begin{cases}
               1, & \mbox{if } i=j~and~i, j\in\{1, 2, 3, 4, 5\} \\
               0, & \mbox{otherwise}.
             \end{cases}
  \]
  Let $\eta$ be a 1-form defined by $\eta(X)=g(X,e_5)$, for arbitrary $X\in\chi(M)$. Let us define (1,1)-tensor field $\phi$ as:
  \begin{align*}
    \phi(e_1) &= e_3, & \phi(e_2) &= e_4, & \phi(e_3) &= -e_1, & \phi(e_4) &= -e_2, & \phi(e_5) &= 0.
  \end{align*}
    Then it satisfy the relations $\eta(\xi)=1$, $\phi^2(X)=-X+\eta(X)\xi$ and $g(\phi X,\phi Y)=g(X,Y)-\eta(X)\eta(Y)$, where $\xi=e_5$ and $X,Y$ is arbitrary vector field on $M$. So, $(M,\phi,\xi,\eta,g)$ defines an almost contact structure on $M$.\par
     We can now deduce that,
    \begin{align*}
      [e_1, e_2] &=0 & [e_1, e_3] &=0 & [e_1, e_4] &=0 & [e_1, e_5] &=e_1  \\
      [e_2, e_1] &=0 & [e_2, e_3] &=0 & [e_2, e_4] &=0 & [e_2, e_5] &=e_2  \\
      [e_3, e_1] &=0 & [e_3, e_2] &=0 & [e_3, e_4] &=0 & [e_3, e_5] &=e_3  \\
      [e_4, e_1] &=0 & [e_4, e_2] &=0 & [e_4, e_3] &=0 & [e_4, e_5] &=e_4  \\
      [e_5, e_1] &=-e_1 & [e_5, e_2] &=-e_2 & [e_5, e_3] &=-e_3 & [e_5, e_4] &=-e_4.
    \end{align*}
Let $\nabla$ be the Levi-Civita connection of $g$. Then from $Koszul's formula$ for arbitrary $X, Y, Z\in\chi(M)$ given by:\\
\begin{eqnarray}
2g(\nabla_XY, Z)&=&Xg(Y, Z)+Yg(Z, X)-Zg(X, Y)-g(X, [Y, Z])\nonumber\\
&-&g(Y,[X, Z])+g(Z, [X, Y]),\nonumber
\end{eqnarray}
we can have:
\begin{align*}
  \nabla_{e_1}e_1 &=-e_5 & \nabla_{e_1}e_2 &=0 & \nabla_{e_1}e_3 &=0 & \nabla_{e_1}e_4 &=0 & \nabla_{e_1}e_5 &=e_1 \\
  \nabla_{e_2}e_1 &=0 & \nabla_{e_2}e_2 &=-e_5 & \nabla_{e_2}e_3 &=0 & \nabla_{e_2}e_4 &=0 & \nabla_{e_2}e_5 &=e_2 \\
  \nabla_{e_3}e_1 &=0 & \nabla_{e_3}e_2 &=0 & \nabla_{e_3}e_3 &=-e_5 & \nabla_{e_3}e_4 &=0 & \nabla_{e_3}e_5 &=e_3 \\
  \nabla_{e_4}e_1 &=0 & \nabla_{e_4}e_2 &=0 & \nabla_{e_4}e_3 &=0 & \nabla_{e_4}e_4 &=-e_5 & \nabla_{e_4}e_5 &=e_4 \\
  \nabla_{e_5}e_1 &=0 & \nabla_{e_5}e_2 &=0 & \nabla_{e_5}e_3 &=0 & \nabla_{e_5}e_4 &=0 & \nabla_{e_5}e_5 &=0.
\end{align*}
Therefore $(\nabla_X\phi)Y=g(\phi X, Y)\xi-\eta(Y)\phi X$ is satisfied for arbitrary $X, Y\in\chi(M)$. So $(M,\phi, \xi, \eta, g)$ becomes a Kenmotsu manifold.\\
The non-vanishing components of curvature tensor are:
\begin{align*}
  R(e_1, e_2)e_2&=-e_1 & R(e_1, e_3)e_3&=-e_1 & R(e_1, e_4)e_4&=-e_1 \\
  R(e_1, e_5)e_5&=-e_1 & R(e_1, e_2)e_1&=e_2 & R(e_1, e_3)e_1&=e_3 \\
  R(e_1, e_4)e_1&=e_4 & R(e_1, e_5)e_1&=e_5 & R(e_2, e_3)e_2&=e_3 \\
  R(e_2, e_4)e_2&=e_4 & R(e_2, e_5)e_2&=e_5 & R(e_2, e_3)e_3&=-e_2 \\
  R(e_2, e_4)e_4&=-e_2 & R(e_2, e_5)e_5&=-e_2 & R(e_3, e_4)e_3&=e_4 \\
  R(e_3, e_5)e_3&=e_5 & R(e_3, e_4)e_4&=-e_3 & R(e_4, e_5)e_4&=e_5 \\ R(e_5, e_3)e_5&=e_3 & R(e_5, e_4)e_5&=e_4.
\end{align*}
Now from the above results we have, $S(e_i,e_i)=-4$ for $i=1,2,3,4,5$  and
\begin{equation}\label{ex1 1}
  S(X, Y)=-4g(X,Y)~~\forall X,Y\in\chi(M).
\end{equation}
Contracting this we have $r=\sum_{i=1}^{5}S(e_i, e_i)=-20=-2n(2n+1)$ where dimension of the manifold $2n+1=5$. Also, we have
\[
  S^*(e_i, e_i)=\begin{cases}
               -1, & \mbox{if } i=1, 2, 3, 4 \\
               0, & \mbox{if } i=5.
               \end{cases}
\]
So,
\begin{equation}\label{ex1 2}
  S^*(X, Y)=-g(X, Y)+\eta(X)\eta(Y)~~\forall X, Y\in\chi(M).
\end{equation}
Hence,
\begin{equation}\label{7.4}
  r^* =\Tr(S^*)=-4.
\end{equation}
Now, we consider a vector field $V$ as
\begin{equation}\label{ex1 3}
  V=x\frac{\partial}{\partial x}+y\frac{\partial}{\partial y}+z\frac{\partial}{\partial z}+u\frac{\partial}{\partial u}+\frac{\partial}{\partial v}.
\end{equation}
Then from the above results we can justify that
\begin{equation}\label{ex1 4}
  (\mathcal{L}_Vg)(X, Y)=4\{g(X,Y)-\eta(X)\eta(Y)\},
\end{equation}
which holds for all $X,Y\in\chi(M)$. Hence from \eqref{ex1 4}, we have
\begin{equation}\label{7.5}
  \sum_{i=1}^{5} (\pounds_V g) (e_i,e_i)=16.
\end{equation}
Now putting $X = Y =e_i$ in the \eqref{1.9}, summing over $i=1, 2, 3, 4, 5$ and using \eqref{7.4} and \eqref{7.5}, we obtain,
\begin{equation}\label{7.6}
  \Lambda=\frac{4\alpha-10\beta-8k}{5}.
\end{equation}
As this $\Lambda$, defined as above satisfies \eqref{3.6}, so $g$ defines a $*$-$k$-Ricci-Yamabe soliton on the 5-dimensional Kenmotsu manifold $M$.\\
Also we can state,
\begin{rem}
\textbf{Case-I:} When $\alpha=1, \beta = 0$, \eqref{7.6} gives $\Lambda=\frac{4-8k}{5}$ and hence $(g, V, \Lambda)$ is a $*$-$k$-Ricci soliton which is shrinking when $k>\frac{1}{2}$, expanding when $k<\frac{1}{2}$ and steady if $k=\frac{1}{2}$.\\\\
\textbf{Case-II:} When $\alpha=0, \beta = 2$, \eqref{7.6} gives $\Lambda=\frac{-20-8k}{5}$ and hence $(g, V, \Lambda)$ is a $*$-$k$-Yamabe soliton which is shrinking if $k>-\frac{5}{2}$, expanding if $k<-\frac{5}{2}$ and shrinking if $k=-\frac{5}{2}$.\\\\
\textbf{Case-III:} When $\alpha=1, \beta = 1$, \eqref{7.6} gives $\Lambda=\frac{-6-8k}{5}$ and hence $(g,V,\Lambda)$ is a $*$-$k$-Einstein soliton which is shrinking if $k>-\frac{3}{4}$, expanding if $k<-\frac{3}{4}$ and steady if $k=-\frac{3}{4}$.
\end{rem}

\end{document}